\documentstyle[12pt]{article}

\newcommand\eq{\begin{equation}}
\newcommand\en{\end{equation}}
\newcommand\eqs{\begin{eqnarray*}}
\newcommand\ens{\end{eqnarray*}}
\def\numberlikeadb{\global\def\theequation
{\thesection.\arabic{equation}}}
\numberlikeadb
\newtheorem{theorem}{Theorem}[section]
\newtheorem{lemma}[theorem]{Lemma}
\newtheorem{corollary}[theorem]{Corollary}
\newtheorem{remark}[theorem]{Remark}
\newcommand\eqa{\begin{eqnarray}}
\newcommand\ena{\end{eqnarray}}
\def\ux{^{(x)}}
\def\xx{{\cal X}}
\def\bone{{\bf 1}}
\def\Blm{\left|}
\def\Brm{\right|}
\def\Bl{\left(}
\def\Br{\right)}
\def\Blb{\left\{}
\def\Brb{\right\}}
\def\Ref#1{(\ref{#1})}
\def\law{{\cal L}}
\def\real{\text{\rm I\kern-2pt R}}
\def\re{\real}
\def\expec{\text{\rm I\kern-2pt E}}
\def\ex{\expec}
\def\prob{\text{\rm I\kern-2pt P}}
\def\pr{\prob}
\def\l{\lambda}
\def\L{\Lambda}
\def\a{\alpha}
\def\b{\beta}
\def\g{\gamma}

\def\e{\varepsilon}
\def\h{\eta}
\def\r{\rho}
\def\NB{{\rm NB\,}}
\def\nin{\noindent}
\def\bsk{\bigskip}
\def\msk{\medskip}
\def\integ{{\bf Z}}
\def\un{^{(n)}}
\def\uj{^{(j)}}
\def\un{^{(n)}}
\def\Po{\text{\rm Po\,}}

\def\text{}
\def\ignore#1{}
\def\th{\theta}
\def\non{\nonumber}
\def\proof{\bsk\nin{\bf Proof.\ }}
\def\epr{\hfil$\diamond$\bsk\bsk}

\def\sjbn{\sum_{j=b+1}^n}
\def\sjbln{\sum_{j=b+1}^{l\wedge n}}
\def\sjblnt{\sum_{j=b+1}^{\lfloor (l/2)\wedge n \rfloor}}
\def\sjbl{\sum_{j=b+1}^l}
\def\sjlt{\sum_{j=1}^\lot}
\def\sjlf{\sum_{j=1}^{\lfloor l/4 \rfloor}}
\def\sjl{\sum_{j=1}^l}
\def\lot{{\lfloor l/2 \rfloor}}
\def\sjltl{\sum_{j=\lot + 1}^l}
\def\sjn{\sum_{j=1}^n}
\def\sji{\sum_{j\ge1}}
\begin{document}

\title{Random combinatorial structures: the convergent case}
\author{A. D. Barbour\footnote{Angewandte Mathematik,
Winterthurerstrasse 190, 8057 Z\"urich, Switzerland;\hfil\break
e-mail:adb@amath.unizh.ch.}\ ~and
Boris L. Granovsky\footnote{Department of Mathematics, Technion, Haifa,
Israel;\hfil\break e-mail: mar18aa@techunix.technion.ac.il.}
\\
Universit\"at Z\"urich and Technion, Haifa}

\maketitle

\begin{abstract}
This paper studies the distribution of the component spectrum of
combinatorial structures such as uniform random forests, in which
the classical generating function for the numbers of (irreducible)
elements of the different sizes converges at the radius of convergence;
here, this property is expressed in terms of the expectations of
{\it independent\/} random variables~$Z_j$, $j\ge1$, whose joint
distribution, conditional on the event that $\sum_{j=1}^n jZ_j = n$,
gives the distribution of the component spectrum for a random structure
of size~$n$.  For a large class of such structures, we show that the
component spectrum is asymptotically composed of~$Z_j$ components
of small sizes~$j$, $j\ge1$, with the remaining part, of size
close to~$n$, being made up of a single, giant component.
\end{abstract}

\section{Introduction}
In this paper, we consider the distribution of the asymptotic
component spectrum of certain decomposable
random combinatorial structures. A structure of size~$n$ is
composed of parts whose (integer) sizes sum to~$n$; we let
$C\un := (C_1\un,C_2\un,\ldots,C_n\un)$ denote its component
spectrum, the numbers of components
of sizes $1,2,\ldots,n$, noting that we always have $\sjn jC_j\un = n$.
For each given~$n$, we assume that the probability distribution
on the space of all
such component spectra satisfies the Conditioning Relation:
\eq\label{CR}
\law(C\un) = \law\Bl (Z_1,Z_2,\ldots,Z_n)\,\bigg|\, \sjn jZ_j = n \Br,
\en
where $Z := (Z_j,\,j\ge1)$ is a sequence of independent random variables,
the same for all~$n$; that is, for $y_1,y_2,\ldots,y_n \in \integ_+$,
\eqa
\lefteqn{\pr[(C_1\un,C_2\un,\ldots,C_n\un) = (y_1,y_2,\ldots,y_n)]}\non\\
&&  = \Blb\pr\left[\sum_{j=1}^n jZ_j=n\right]\Brb^{-1}
  \prod_{j=1}^n \pr[Z_j = y_j]\,\bone_{\{n\}}
  \Bl \sum_{j=1}^n jy_j \Br.\phantom{HHH} \label{CR2}
\ena
This apparently curious
assumption is satisfied by an enormous number of classical
combinatorial objects, such as, for instance, permutations of~$n$
objects under the uniform distribution, decomposed into cycles
as components,
when the~$Z_j$ are Poisson distributed, with $Z_j \sim \Po(1/j)$;
or forests of unlabelled unrooted trees under the uniform
distribution, decomposed into tree components, when the~$Z_j$
are negative binomially distributed: see Arratia, Barbour
\& Tavar\'e~(2002, Chapter~2)[ABT]
for many more examples.  However, such structures
also arise in other contexts.  For instance, the state of a
coagulation--fragmentation process evolving in a collection of~$n$
particles can be described by the numbers~$C_j\un$ of clusters
of size~$j$, $1\le j\le n$, and if such a process is reversible
and Markov, then its equilibrium distribution satisfies the
conditioning relation for some sequence~$Z$ of random variables.
In particular, under mass action kinetics,
it follows that~$Z_j\sim \Po(a_j)$, where $(a_j,\,j\ge1)$ are
positive reals, determined by the coagulation
and fragmentation rates; see Whittle~(1965), Kelly~(1979, Chapter~8),
Durrett, Granovsky \& Gueron~(1999)  and Freiman \& Granovsky~(2002a).

In order to describe the asymptotics as $n\to\infty$, it
is necessary first to say something about how the distributions of
the~$Z_j$ vary with~$j$.  Now the distribution given in~\Ref{CR2}
remains the same
if the random variables~$Z_j$ are replaced by `tilted'
random variables~$Z_j\ux$, where
\eq\label{tilt-Z}
\pr[Z_j\ux = i] = \pr[Z_j = i] x^{ji} / k_j(x),
\en
for {\it any}~$x>0$ such that
$$
k_j(x) := \ex\Blb x^{jZ_j} \Brb < \infty.
$$
Specializing to the setting in which~$Z_j \sim\Po(a_j)$ for each~$j$,
this means that exactly the same distributions are obtained for
each~$n$ in~\Ref{CR2} if~$a_j$ is replaced by~$a_jx^j$ for each~$j$,
for any fixed~$x>0$. Thus geometrically fast growth or decay of
the~$a_j$ can be offset by choosing $x^{-1} = \lim_{j\to\infty}
a_j^{1/j}$ (should the limit exist), without changing the asymptotics.
Hence, to find an interesting range of possibilities, we look at
rates of growth or decay of~$\ex Z_j$ which (if necessary, after
appropriate tilting) can be described by a power law:
$\ex Z_j \sim Aj^{\a}$ as $j\to\infty$, or, more generally,
$\ex Z_j$ regularly varying with exponent~$\a \in \re$.

Three ranges of~$\a$ can then broadly be distinguished.  The most
intensively studied is that where $\a=-1$, and within this
the logarithmic class, in which $\ex Z_j \sim \pr[Z_j=1]
\sim \th j^{-1}$, for some~$\th > 0$:  see the book [ABT] for
a detailed discussion.  For $\a > -1$, the expansive case,
the asymptotics were explored for Poisson distributed~$Z_j$
in Freiman \& Granovsky~(2002a,b),
with the help of Khinchine's probabilistic
method, and particular models have been studied by many authors.
Here, we treat the convergent case,
in which $\a < -1$, in considerable generality.  Our approach is
quite different from the classical approach by way of generating
functions, thereby allowing distributions other than the standard
Poisson and negative binomial to be easily discussed.
Note also that not all classical combinatorial structures fall
into one of these three categories: random set partitions,
studied using the Conditioning Relation by Pittel~(1997),
have Poisson distributed~$Z_j$ with means $x^j/j!$, which
are never regularly varying, whatever the choice of~$x>0$.

As will be seen in
what follows, a key element in the arguments is establishing
the asymptotics of the probabilities $\pr[T_{bn}(Z) = l]$
for~$l$ near~$n$, where, for $y := (y_1,y_2,\ldots)\in \integ_+^\infty$,
\eq\label{Tbn-def}
T_{bn}(y) := \sjbn jy_j,\qquad 0 \le b < n.
\en
That this should be so is clear from~\Ref{CR2},
in which the normalizing constant is just the probability
$\pr[T_{0n}(Z)=n]$,
and is the only element which cannot immediately be written down.
In the context of reversible coagulation--fragmentation processes
with mass--action kinetics, the partition
function~$c_n$ investigated by Freiman \& Granovsky~(2002a) is
given by
\eq\label{partition}
c_n := \exp\Blb \sjn a_j \Brb \pr[T_{0n}(Z)=n],
\en
explaining its relation to many of their quantities of interest.
Now, in the expansive case, taking Poisson distributed~$Z_j$ with
means~$a_j \sim Aj^\a$, $\a > -1$, one has
$$
\ex T_{0n}(Z) \asymp n^{2+\a} \gg n \quad {\rm and}\quad
  {\rm SD\,}(T_{0n}(Z)) \asymp n^{(3+\a)/2} \ll \ex T_{0n}(Z).
$$
The Bernstein inequality then implies that, for large~$n$, the probability
$\pr[T_{0n}(Z)=n]$ is extremely small, making a direct asymptotic
argument very delicate.  However, recall from~\Ref{tilt-Z}
that the conditioning
relation~\Ref{CR} delivers the same distribution for the combinatorial
structure if the Poisson distributed random variables~$Z_j$
with means~$a_j$ are replaced by Poisson distributed random variables
$Z_j^{(x)}$ with means $a_j x^j$, for {\it any\/} $x > 0$.
Choosing $x=x_n$ in such a way that $\ex T_{0n}(Z^{(x)}) = n$
makes the probability $\pr[T_{0n}(Z^{(x)}) =n]$ much larger, and
a local limit theorem based on the normal approximation can then
be used to determine its asymptotics.  The resulting
component spectra typically have almost all their weight in components
of size about $n^{1/(\a + 2)}$, a few smaller components making
up the rest.

For the logarithmic case, taking Poisson distributed~$Z_j$ with
means~$a_j \sim \th/j$, $\th > 0$, one has
$$
\ex T_{0n}(Z) \sim n\th \quad {\rm and}\quad {\rm SD\,}(T_{0n}(Z))
  \asymp n,
$$
so that no tilting is required.  However, since
$T_{0n}(Z) \ge 0$,
these asymptotics also imply that $\law(n^{-1}T_{0n}(Z))$ is not close to a
normal distribution --- there is a different limiting distribution that has
a density related to the Dickman function from number theory ---
and special techniques have to be developed in order to complete
the analysis.  Here, the component spectra typically have components
of sizes around~$n^\b$ for {\it all\/} $0\le\b\le1$.

In the convergent case, taking Poisson distributed~$Z_j$ with
means~$a_j \sim Aj^\a$, $\a < -1$, the sequence of random
variables~$T_{0n}(Z)$ converges without normalization, and both the
methods of proof and the typical spectra as $n\to\infty$
are again qualitatively different.  We demonstrate that,
for large~$n$,
the typical picture is that of small components whose numbers
have the independent joint distribution of the~$Z_j$, the remaining
weight being made up by a {\it single\/} component of size close
to~$n$. This remains true without the
Poisson assumption, under fairly weak conditions; for instance,
our theory applies to the
example of uniform random forests, where
the asymptotic distribution of the size of the largest component
was derived using generating function methods by
Mutafchiev~(1998). Bell, Bender, Cameron and Richmond~(2000, Theorem~2)
have also used generating function methods to examine the convergent case
for labelled and unlabelled structures, which, in our setting,
correspond to Poisson and negative binomially distributed~$Z_j$'s,
respectively; we allow an even wider choice of distributions for
the~$Z_j$.  They use somewhat different conditions, and are
primarily interested in whether or not the probability that the
largest component is of size~$n$ has a limit as $n\to\infty$,
though they also consider the limiting distribution of the number
of components. Under our conditions, these limits always exist.

\section{Results}
\setcounter{equation}{0}
We work in a context in which the random variables~$Z_j$
may be quite general, provided that, for large~$j$,
their distributions are sufficiently close to Poisson. From now on,
we use the notation $a_j := \ex Z_j$, and then write
$a_j = j^{-q-1}\l(j)$ for $q = -\a-1 > 0$ in the convergent
case, where the quantities~$\l(j)$ are required to satisfy
certain conditions given below.

Since now $a_j \to 0$, being close to Poisson
mainly involves assuming that
$\pr[Z_j \ge 2] \ll \pr[Z_j = 1]$ as $j\to\infty$, so that the~$Z_j$
can be thought of as independent random variables which usually
take the value~$0$, and occasionally (but only a.s.\ finitely often)
the value~$1$.  This setting is broad enough to include a number
of well known examples, including uniform random forests
consisting of (un)labelled (un)rooted trees.
In such circumstances, we are able to use a technique
based on recurrence relations which are exactly true for Poisson
distributed~$Z_j$, and which can be simply derived using Stein's
method for the compound Poisson distribution (Barbour, Chen and Loh, 1992).
A corresponding
approach is used in [ABT], though the detail of the argument here
is very different.

In describing the closeness of the distributions of the~$Z_j$ to Poisson,
we start by exploiting any divisibility that they may possess,
supposing that each~$Z_j$ can be
written in the form $Z_j = \sum_{k=1}^{r_j} Z_{jk}$ for some
$r_j\ge1$, where, for each~$j$, the
non-negative integer valued random variables
$(Z_{jk},\,1\le k\le r_j)$ are independent and identically distributed.
Clearly, this is always possible if we take $r_j=1$. However, Poisson
distributions are infinitely divisible ($r_j$ may be taken to be arbitrarily
large), and the error bounds in our approximations
become correspondingly smaller, if we are able to choose larger~$r_j$.
Note, however, that negative binomially distributed~$Z_j$
also have infinitely divisible distributions, so that closeness to
Poisson is not a consequence of infinite divisibility alone.
We now define $(\e_{js},\,s,j\ge1)$ by setting
\eqa
r_j\pr[Z_{j1} = 1] &=:& j^{-q-1}\l(j) (1 - \e_{j1});\non\\
r_j\pr[Z_{j1} = s] &=:& j^{-q-1}\l(j)\e_{js},\quad s\ge2,
  \label{Z-dist}
\ena
so that then
$$
0 \le \e_{j1} = \sum_{s\ge2} s\e_{js} \le 1,
$$
because $j^{-q-1}\l(j) = a_j = \ex Z_j = r_j\ex Z_{j1}$. We then assume that
\eq\label{epsilon-A1}
0 \le \e_{js} \le \e(j)\g_s, \quad s\ge2,
\en
where
\eq\label{epsilon-A2}
G  := \sum_{s\ge2}s\g_s < \infty \quad {\rm and}
\quad \lim_{j\to\infty}\e(j) = 0;
\en
we write $\e^*(j) := \max_{l\ge j+1} \e(l)$ and $r^*(j) :=
\min_{l > j} r_l$.
For the subsequent argument, we need to strengthen
\Ref{epsilon-A2} by assuming in addition that
\eq\label{epsilon-A3}
G_q := \sum_{s\ge2}L_s s^{1+q}\g_s < \infty,\quad {\rm where} \quad
L_s := \sup_{l\ge s}\{\l(\lfloor l/s \rfloor)/\l(l)\}.
\en

We also need some conditions on the function~$\l$. We assume that
\eqa
\l^+(l) &:=& \max_{1\le s\le l}\l(s) = o(l^\b) \quad
               {\rm for\ any} \ \b > 0;\label{lambda-A1}\\
L &:=& \sup_{l\ge2}\max_{l/2 < t \le l}\{\l(l-t)/\l(l)\} < \infty,
   \label{lambda-A2}
\ena
and that
\eqa
&&\lim_{l\to\infty}\{\l(l-s)/\l(l)\} = 1 \quad {\rm for\ all}\ s \ge 1;
    \label{lambda-A3}
\ena
note that, if~$\l$ is slowly varying at infinity,
then conditions \Ref{lambda-A1}--\Ref{lambda-A3} are
automatically satisfied, and that~$L_s$ defined in~\Ref{epsilon-A3}
is finite.
We then write $\L_\b := \max_{l\ge1} l^{-\b}\l(l)$ for $\b > 0$,
and we also observe that
\eq\label{T0inf-finite}
\pr[Z_{jk} \ge 1 \ \mbox{for}\ \infty\ \mbox{many}\ j,k] = 0 \quad
 \mbox{and hence that}\quad T_{0\infty}(Z) < \infty \ \mbox{a.s.},
\en
from \Ref{Z-dist}, \Ref{lambda-A1} and the Borel--Cantelli lemma.
  Finally,
we assume that the distributions of the random variables~$Z_{j1}$
of~\Ref{Z-dist} are such that
\eq\label{P0-bnd}
p_0 := \min_{j\ge1}\pr[Z_{j1}=0] > 0.
\en
This restriction can actually be dispensed with --- see Remark~\ref{p0-rk}
--- but it makes the proofs somewhat simpler.

We are now in a position to state our first theorem, in which the
asymptotics of the probabilities  $\pr[T_{bn}(Z) = l]$ are
described.

\begin{theorem}\label{LLT}
Suppose that conditions \Ref{epsilon-A1} --~\Ref{lambda-A3}
are satisfied for some $q>0$ and that~\Ref{P0-bnd} holds.
For $1 \le l \le n$, define
$$
 H_n(l) := \max_{0\le b\le l-1}|\l^{-1}(l)l^{1+q}\pr[T_{bn} = l] - 1|.
$$
Then $H(l) := \sup_{n\ge l}H_n(l)$ satisfies $\lim_{l\to\infty}H(l) = 0$.
\end{theorem}

\nin Note that the condition $G_q < \infty$ of~\Ref{epsilon-A3} is
really needed here: see Remark~\ref{Gq-rk}.

As is strongly suggested by the formula~\Ref{CR2},
Theorem~\ref{LLT}, in giving the asymptotics of $\pr[T_{0n}(Z)=n]$,
can directly be applied to establish the asymptotic joint distribution
of the entire component spectrum.  This is given in the following
theorem. For probability distributions on a discrete
set~$\xx$, we define the total variation distance~$d_{TV}$ by
$$
d_{TV}(P,Q) := \sup_{A\subset\xx}|P(A)-Q(A)|.
$$

\begin{theorem}\label{joint-dns}
Suppose that conditions \Ref{epsilon-A1} --~\Ref{lambda-A3}
are satisfied for some $q>0$, and that~\Ref{P0-bnd}
holds.  Then
$$
\lim_{n\to\infty} d_{TV}(\law(C\un),Q_n) \to 0,
$$
where $Q_n$ is the distribution of $(Z_1,Z_2,\ldots,Z_n) + e(n-T_{0n}(Z))$,
and $e(j)$ denotes the $j$'th unit $n$-vector if $j\ge1$, and the zero
$n$-vector otherwise.
\end{theorem}

Theorem~\ref{joint-dns} has a number of immediate
consequences, which all follow directly because $T_{0\infty}(Z)
< \infty$ a.s.

\begin{corollary}\label{last-cor}\hfil\break
(a) For any fixed $k\ge 1$,
$$
\law(C_1\un,\ldots,C_{k}\un)\to \law(Z_1,\ldots, Z_k)
  \quad \mbox{as} \quad  n\to \infty.
$$

\nin(b) If $Y_n:=\max\{j:C_j\un > 0\}$ and
$K_n:=\min\{j:C_j\un > 0\}$ are the sizes of the maximal
and minimal components of the spectrum, then, as $n\to\infty$,
$$
\law(n-Y_n) \to \law(T_{0\infty}(Z))
$$
and
$$
\pr[K_n > b] \to \prod_{j=1}^b \pr[Z_j=0]
$$
for any $b > 1$. In particular, it follows that
\eqa
\lim_{n\to \infty}\pr[Y_n = K_n =n]=\prod_{j\ge1}\pr[Z_j=0].\label{1}
\ena

\nin(c) The asymptotic distribution of the number of
components~$X_n$ of the spectrum is given by
$$
\law(X_n) \to \law\Bl 1 + \sum_{j\ge1} Z_j \Br.
$$
\end{corollary}

\begin{remark}\label{final}
\textnormal{
The assertion~(a) of the above corollary states the
asymptotic independence of the numbers of components of small
sizes, a fact that has also been established in~[ABT] in the logarithmic
case, and also in the Poisson setting for
$q>0$ in Freiman and Granovsky~(2002b).
This fact can be viewed as a particular manifestation of the
heuristic general principle of asymptotic
independence of particles in models of statistical physics.}

\textnormal{
Assertion (b) says that, as $n\to \infty$,  the structures
considered exhibit the gelation phenomenon; the formation, with
positive probability, of a component with size comparable to~$n$
 (see, for example, Whittle~(1986), Ch.~13).
Gelation also occurs in the logarithmic case [ABT], while it
is not seen for $q>0$ in the setting of Freiman and Granovsky~(2002b).
In this sense, $q=0$ ($\a=-1$) represents a critical value
of the exponent.}
\end{remark}

\ignore{
First, if~$Y_n$ denotes the size of the maximal, `giant'
component, it follows that
$$
\law(n-Y_n) \to \law(T_{0\infty}).
$$
Similarly, if $Q_n$ denotes the size of the smallest component, then
$$
\pr[Q_n > b] \to \prod_{j=1}^b \pr[Z_j=0].
$$
The two are equal, and equal to~$n$, with asymptotic probability
$$
\pr[T_{0\infty}=0] = \pr[Y_n=n] = \prod_{j\ge1} \pr[Z_j=0].
$$
The asymptotic distribution of the number of components~$X_n$ is
given by $\law\Bl 1 + \sum_{j\ge1} Z_j \Br$.
}

\bsk
Now $\pr[Y_n = n]$
is the probability that a structure is `connected', as, for instance,
in Bell, Bender, Cameron and Richmond~(2000), who give a very
general discussion
of circumstances in which $\rho := \lim_{n\to\infty}\pr[Y_n=n]$ exists,
as well as giving a formula for the asymptotic distribution of~$X_n$.
They work in the settings of either labelled or unlabelled structures; in
our terms, they assume that the~$Z_j$ have either Poisson or negative
binomial distributions, respectively. Theorem~\ref{joint-dns} implies that
$\rho$ always exists under our conditions, and gives its value.

\bsk

\nin{\bf Example.}\ We apply our results to some classical models of
random forests, referring for a discussion of the literature to
the books of Pavlov~(2000) and Kolchin~(1998); see also
Mutafchiev~(1998, pp.~212--213).  We begin by
considering the uniform distribution over
all forests of unlabelled, unrooted trees.  The number~$m_j$ of
such trees of size~$j$ was studied by Otter~(1948), who showed
that $m_j \sim c\r^{-j}j^{-5/2}$, where $\r < 1$, and gave values
for both $\r$ and~$c$.  This combinatorial structure satisfies
the conditioning relation with negative binomial
random variables $Z_j \sim \NB(m_j,\r^j)$, so that
$$
\pr[Z_j = s] = (1-\r^j)^{m_j}{m_j+s-1 \choose s}\r^{js},\quad s\ge0.
$$
It thus follows that
$\ex Z_j = m_j \r^j/(1-\r^j) \sim cj^{-5/2}$, implying that our
results can be applied with $\lambda(j)\to c$ and $q=3/2$.
Note that, if we take $r_j=1$ for all~$j$, we have
$$
\pr[Z_j=2] = (1-\r^j)^{m_j}{m_j+1\choose2}\r^{2j}
   \asymp (m_j\r^j)^2,
$$
so that $\e_{j2} \asymp j^{-5/2}$ as $j\to\infty$.  On the other
hand, negative binomial distributions are infinitely divisible,
and other choices of~$r_j$ in~\Ref{Z-dist} are possible: for each~$j$,
we can take
$Z_{jk} \sim \NB(m_j/r_j,\r^j)$, $1\le k\le r_j$, for any choice of~$r_j$.
The corresponding values of~$\e_{js}$, $s\ge2$, are then given,
using~\Ref{Z-dist}, by
\eqs
r_j\pr[Z_{j1}=s] &=&
  r_j(1-\r^j)^{\frac{m_j}{r_j}}{\frac{m_j}{r_j}+s-1\choose s}\r^{sj}\\
&=&r_j(1-\r^j)^{\frac{m_j}{r_j}}\r^{sj}
  \frac{(\frac{m_j}{r_j}+s-1)\cdots (\frac{m_j}{r_j}+1)\frac{m_j}{r_j}}{s!}\\
&=&\{m_j \r^j/(1-\r^j)\}\epsilon_{js},
\ens
from which, for fixed~$j$ and $s\ge2$, we deduce the limiting value
$$
\epsilon_{js}^* = s^{-1}(1-\r^j)\r^{(s-1)j}
$$
of $\epsilon _{js}$ as $r_j\to \infty$.
Note that, as $j\to\infty$,
$\e_{j2}^* \sim 2^{-1}\r^{j}$ is of very much smaller order
than the order~$j^{-5/2}$ obtained for~$\e_{j2}$ when taking $r_j=1$.
As a result, many of the contributions to the bound~$H(l)$
of Theorem~\ref{LLT} for the relative error
in approximating $\pr[T_{bn}=l]$ are reduced. These include the
terms arising from
$\h_0'$, $\h_1$ and~$\h_2$, which enter in \Ref{last-line}
and~\Ref{2nd-line-term1} below; furthermore, as observed in
Remark~\ref{id-rk}, letting $r_j \to \infty$ also allows us to
take $p_0=1$ and $\h_k(l) = 0$, $3\le k\le 6$.

Similar arguments can be used for forests of unlabelled, rooted
trees, now with $m_j \sim c'\r^{-j}j^{-3/2}$.
For forests of labelled, (un)rooted trees, $\law(T_{0\infty})$ is
the compound Poisson distribution of $\sji jZ_j$, where
$$
Z_j \sim \Po\Bl {j^{j-2}\over j!e^j} \Br \quad\mbox{(unrooted)};
\qquad
Z_j \sim \Po\Bl {j^{j-1}\over j!e^j} \Br \quad\mbox{(rooted)}.
$$
The asymptotics of $\law(n-Y_n)$ then implied by Corollary~\ref{last-cor}
do not appear to agree with those of Mutafchiev~(1998).

\section{Proofs}
\setcounter{equation}{0}

\subsection{The perturbed Stein recursion and the basic lemma}

Stein's method for the Poisson distribution~$\Po(a)$ is based on
the Stein--Chen identity
$$
\ex\{Zf(Z)\} = a\ex f(Z+1),
$$
true for all bounded functions $f:\integ_+ \to \re$ when $Z\sim\Po(a)$;
this can be checked by writing the expectations on each
side of the equation as sums, and then
examining the coefficients of~$f(l)$  for each~$l\ge0$.
In particular, it then follows that
$$
\ex\{jZ_jg(jZ_j)\} = ja_j\ex g(jZ_j + j)
$$
if $Z_j \sim \Po(a_j)$, by putting $f(l) = g(jl)$.  Hence, for the
compound Poisson distributed weighted sum
$$
T^*_{bn} := T_{bn}(Z) = \sjbn jZ_j,
$$
when $Z_j \sim \Po(a_j)$ and the~$Z_j$ are all independent, we
deduce the Stein identity
\eq\label{Stein-identity}
\ex \{T^*_{bn}g(T^*_{bn})\} = \sjbn ja_j \ex g(T^*_{bn}+j),
\en
true for all bounded functions $g:\integ_+ \to \re$
and for any $0\le b < n$.
Taking $g = \bone_{\{l\}}$, for any $l \ge b+1$, it thus follows that
\eqa
l\pr[T^*_{bn}=l] &=& \sjbn j^{-q}\l(j)\pr[T^*_{bn} = l-j] \non \\
  &=& \sjbln j^{-q}\l(j)\pr[T^*_{bn} = l-j], \quad l \ge b+1;
     \label{Poisson-probs}
\ena
note that this recursion can also be deduced directly by differentiating
the compound Poisson generating function, and equating coefficients.
Recursion~\Ref{Poisson-probs}, coupled with the fact that
$\pr[T^*_{bn} = l] = 0$
for $1\le l \le b$, successively expresses the probabilities
$\pr[T^*_{bn}=l]$ in terms of the probability
$\pr[T^*_{bn}=0]$.  In particular,
if~$l\le n$ is large and if $\{j^{-q}\l(j)\}/\{l^{-q}\l(l)\}$ is close
to~$1$ when~$j$ is close to~$l$, it suggests that
$$
l\pr[T^*_{bn}=l] \approx l^{-q}\l(l) \pr[T^*_{bn} < l-b-1] \approx
  l^{-q}\l(l),
$$
giving the large~$l$ asymptotics for $\pr[T^*_{bn}=l]$.  Our
approach consists of turning this heuristic into a precise
argument, which can be applied also when the~$Z_j$ do not have
Poisson distributions.

Observing that the
Stein identity~\Ref{Stein-identity} is deduced from the
Stein--Chen identity
\eq\label{Stein-Poisson}
\ex\{Z_j g(T_{bn}^*)\} = j^{-1-q}\l(j) \ex\{g(T_{bn}^*+j)\},
  \quad b+1 \le j \le n,
\en
when $Z_j \sim \Po(j^{-1-q}\l(j))$,
our first requirement is to establish an analogue of
\Ref{Stein-Poisson} for more general random variables~$Z_j$.
To do so, as in the previous section,
we suppose that each~$Z_j$ can be
written in the form $Z_j = \sum_{k=1}^{r_j} Z_{jk}$ for some
$r_j\ge1$, where, for each~$j$, the
non-negative integer valued random variables
$(Z_{jk},\,1\le k\le r_j)$ are independent and identically distributed.
Then, writing $T_{bn} := T_{bn}(Z)$, it is immediate that
$$
\ex\{Z_{j1} g(T_{bn})\}
  = \sum_{s\ge1} s\pr[Z_{j1} = s]\ex g(T_{bn}\uj + js),
$$
where $T_{bn}\uj := T_{bn} - jZ_{j1}$, so that, with the
above definitions,
\eqa
\ex\{T_{bn} g(T_{bn})\} &=& \sjbn \ex\{jZ_jg(T_{bn})\}\label{star}\\
&=& \sjbn jr_j
    \sum_{s\ge1} s\pr[Z_{j1} = s]\ex g(T_{bn}\uj + js) \non\\
&=& \sjbn j^{-q}\l(j) \ex g(T_{bn} + j) \non\\
&&\quad\mbox{} +
   \sjbn j^{-q}\l(j) \{(1-\e_{j1})\ex g(T_{bn}\uj + j)
                                    - \ex g(T_{bn} + j)\}\non\\
&&\quad\mbox{} +
    \sjbn \sum_{s\ge2} j^{-q}\l(j)s\e_{js}\ex g(T_{bn}\uj + js).
   \label{Stein-identity-2}
\ena
Taking $g = \bone_{\{l\}}$ as before then gives the recursion
\eqa
l\pr[T_{bn}=l] &=& \sjbln j^{-q}\l(j)\pr[T_{bn} = l-j]\non\\
&&\ \mbox{} +
   \sjbln j^{-q}\l(j) \{(1-\e_{j1})\pr[T_{bn}\uj =l-j]
                                    - \pr[T_{bn} = l-j]\}\non\\
&&\ \mbox{} +
    \sjblnt \sum_{s\ge2} j^{-q}\l(j)s\e_{js}\pr[T_{bn}\uj = l - js],
   \label{general-probs}
\ena
which can be understood as a perturbed form of the
recursion~\Ref{Poisson-probs}.

In order to show that the perturbation
is indeed small, it is first necessary to derive bounds for the
probabilities $\pr[T_{bn} = s]$ and $\pr[T_{bn}\uj = s]$.
However, since $\pr[T_{bn}=s] \ge \pr[Z_{j1}=0]\,\pr[T_{bn}\uj = s]$,
we have the immediate bound
\eq\label{Tbnj-bnd}
\pr[T_{bn}\uj = s] \le p_0^{-1}\pr[T_{bn} = s],\quad s=0,1,\ldots,
\en
where $p_0 > 0$ is as in~\Ref{P0-bnd}.
Hence the following lemma is all that is required.

\begin{lemma}\label{upper-bnd}
Suppose that conditions \Ref{epsilon-A1} --~\Ref{lambda-A3}
are satisfied for some $q>0$, and that~\Ref{P0-bnd} holds.
Then there exists a
constant~$K>0$, depending only on the distributions of the~$Z_j$, such
that
$$
\pr[T_{bn}=l] \le K\l(l)l^{-1-q},\qquad l\ge1.
$$
\end{lemma}

\proof
For $1\le l\le b$, the statement is trivial.  For larger~$l$, we proceed
by induction, using the recursion~\Ref{general-probs}, in which, on the
right hand side, probabilities of the form $\pr[T_{bn}=s]$ appear only
for $s<l$, so that we may suppose that then
$\pr[T_{bn}=s] \le K\l(s)s^{-1-q}$ for all $1\le s < l$.
Under this hypothesis, we split
the right hand side of~\Ref{general-probs} into three terms,
which we bound separately;
we take the first two lines together, and then split the third
according to the value taken by~$js$.

For the first term, we use~\Ref{Tbnj-bnd}, the induction hypothesis
and conditions \Ref{lambda-A1} and~\Ref{lambda-A2} to give
\eqa
\lefteqn{\sjbln j^{-q}\l(j) (1-\e_{j1})\pr[T_{bn}\uj =l-j]}\non\\
&\le& \sjlt j^{-q}\l(j)p_0^{-1}\pr[T_{bn} =l-j] +
  \sjltl  j^{-q}\l(j)p_0^{-1}\pr[T_{bn} =l-j] \non\\
&\le&  p_0^{-1} \l^+(\lot) KL\l(l)(2/l)^{1+q} \sjlt j^{-q} +
  p_0^{-1} L\l(l)(2/l)^q  \non\\
&=&  K\l(l)l^{-q}\h_0(l) + p_0^{-1} L\l(l)(2/l)^q,\label{term1}
\ena
where
$$
\h_0(l) :=  p_0^{-1} 2^{1+q} \l^+(\lot) L l^{-1} \sjlt j^{-q} = o(1)
\quad{\rm as}\ l\to\infty.
$$
For the second term, arguing much as before, we have
\eqa
\lefteqn{\sjblnt \sum_{s\ge2} \bone_{\{js \le \lot\}}
    j^{-q}\l(j)s\e_{js} \pr[T_{bn}\uj = l - js]}  \non\\
&\le&  \sjlt \sum_{s\ge2}  \bone_{\{js \le \lot\}}
    j^{-q}\l(j)s\e_{js} p_0^{-1} KL\l(l)(2/l)^{1+q} \non\\
&\le&  \l(l)l^{-q}\, p_0^{-1} 2^{1+q} \l^+(\lot) KL l^{-1} \sjlt
    j^{-q}\e(j)G  \non\\
&\le&  \e^*(0) G  K\l(l)l^{-q}\h_0(l). \label{term2}
\ena

For the third and final term, we have
\eqa
\lefteqn{\sjblnt \sum_{s\ge2} \bone_{\{\lot < js \le l\}}
    j^{-q}\l(j)s\e_{js} \pr[T_{bn}\uj = l - js]}  \non\\
&\le& \sum_{s=2}^l
  \sum_{j=\lfloor l/2s \rfloor + 1}^{\lfloor l/s \rfloor - 1}
    j^{-q}\l(j)s\e_{js} \pr[T_{bn}\uj = l - js]
  +  \sum_{s=2}^l \lfloor l/s \rfloor^{-q}\l(\lfloor l/s \rfloor)s
     \e_{\lfloor l/s \rfloor,s}\non\\
&=& S_1 + S_2,
\ena
say.  Now
\eqa
S_1
&\le& \sum_{s=2}^l \sum_{j=\lfloor l/2s \rfloor + 1}^{\lfloor l/s \rfloor - 1}
    j^{-q}\l(j)s\e(j)\g_s p_0^{-1}K\l(l-js) (l-js)^{-1-q}
   \non\\
&\le& p_0^{-1}K \sum_{s=2}^l  (l/2s)^{-q}
       L_sL\l(l)s\g_s\e^*(\lfloor l/2s \rfloor) R_q s^{-1-q/2},
\ena
where $R_q := \L_{q/2}\sum_{t\ge1}t^{-1-q/2}$, and this implies that
\eq\label{term3.1}
S_1 \le K\l(l)l^{-q} \h_1(l),
\en
where
\eqa
\h_1(l) &:=& p_0^{-1}LR_q2^q\,\min_{2 \le t \le l}
   \Blb \e^*(\lfloor l/2t \rfloor) \sum_{s=2}^t s^{q/2}L_s\g_s
      + \e^*(0) \sum_{s \ge t+1} s^{q/2}L_s\g_s\Brb \non\\
 &=& o(1) \quad{\rm as}\ l\to \infty, \non
\ena
in view of \Ref{epsilon-A2} and~\Ref{epsilon-A3}.
For~$S_2$, we have
\eqa
S_2 &\le&
\sum_{s=2}^l \lfloor l/s \rfloor^{-q}\l(\lfloor l/s \rfloor)
     \e(\lfloor l/s \rfloor)s\g_s \non\\
&\le& \l(l)l^{-q} \sum_{s\ge2} s^{1+q}L_s\g_s \e(\lfloor l/s \rfloor)\non\\
&:=& \l(l)l^{-q}\h_2(l),\label{S2-bnd}
\ena
where $\h_2(l) = o(1)$ as $l\to\infty$, again
in view of \Ref{epsilon-A2} and~\Ref{epsilon-A3}.

Collecting these bounds, we can apply \Ref{general-probs} to show that
\eq\label{first-bnd}
l\pr[T_{bn}=l] \le \l(l)l^{-q}\{2^qLp_0^{-1} + \h_2(l)
          + K[\h_0(l)(1+\e^*(0)G )  + \h_1(l)]\};
\en
and this in turn is less than $K\l(l)l^{-q}$ provided that
$$
K\{1 -  [\h_0(l)(1+\e^*(0)G )  + \h_1(l)]\} >  2^qLp_0^{-1} + \h_2(l),
$$
which can be achieved uniformly
for all $l \ge l_0$, for some large~$l_0$, by choosing
$K \ge 2^{q+1}Lp_0^{-1}$. As observed before, $\pr[T_{bn}=l]=0$
for $1\le l\le b$.  For $b+1 \le l \le l_0$, we can suppose that
$\pr[T_{bn} = t] \le K_{l-1} \l(t)t^{-1-q}$ for all $t \le l-1$, and deduce
from~\Ref{first-bnd} that $\pr[T_{bn} = t] \le K_l \l(t)t^{-1-q}$ for all
$t \le l$, if we take
$$
K_l = \max\{K_{l-1},2^qLp_0^{-1} + \h_2(l)
          + K_{l-1}[\h_0(l)(1+\e^*(0)G )  + \h_1(l)]\};
$$
this then completes the proof.  \epr

Lemma~\ref{upper-bnd}, together with the bounds derived in the course of
its proof, are enough to enable us to exploit the
recursion~\Ref{general-probs}, and thereby to prove Theorems \ref{LLT}
and~\ref{joint-dns}; the detailed argument is given in the next two
sections.

\subsection{Proof of Theorem~\ref{LLT}}
We exploit the recursion~\Ref{general-probs},  observing first that the
contribution from its last line
was bounded in the proof of Lemma~\ref{upper-bnd} by
\eq\label{last-line}
 \l(l)l^{-q}\{\e^*(0) G  K\h_0(l) + K\h_1(l) + \h_2(l)\},
\en
uniformly in $0\le b\le l-1$.  We now need to examine the second line
in more detail.  First, note that, by Lemma~\ref{upper-bnd},
for $l\le n$,
\eqa
\lefteqn{\sjbln j^{-q}\l(j) \e_{j1}\pr[T_{bn}\uj =l-j]} \non\\
&\le& \sjlt j^{-q}\l(j)G\e(j) p_0^{-1}(2/l)^{1+q}KL\l(l)
   + p_0^{-1}G\e^*(\lot)L\l(l)(2/l)^q \non\\
&\le& \l(l)l^{-q}(KG\e^*(0)\h_0(l) + \h'_0(l)),\label{2nd-line-term1}
\ena
where
$$
\h'_0(l) := 2^q p_0^{-1} G L \e^*(\lot) = o(1)\quad {\rm as}\ l\to\infty.
$$
The remaining part of the second line of~\Ref{general-probs} is then
bounded by
\eqa
\lefteqn{\Blm \sjbln j^{-q}\l(j) \{\pr[T_{bn}\uj =l-j]
                                    - \pr[T_{bn} = l-j]\}\Brm } \non\\
&=& \Blm \sjbl  j^{-q}\l(j) \{\pr[T_{bn}\uj =l-j]
         - \sum_{s\ge0}\pr[Z_{j1}=s]\pr[T_{bn}\uj = l-j(s+1)]\}\Brm \non\\
&\le& \sjl  j^{-q}\l(j)\{\pr[Z_{j1} \ge 1]\pr[T_{bn}\uj =l-j] \non\\
&&\qquad\mbox{}  + \sum_{s\ge1}\pr[Z_{j1}=s]\pr[T_{bn}\uj =l-j(s+1)]\}.
          \label{2nd-line-term2}
\ena

We now observe, using Lemma~\ref{upper-bnd}, \Ref{Tbnj-bnd},
\Ref{Z-dist} and~\Ref{lambda-A2}, that
\eqa
\lefteqn{\sjl  j^{-q}\l(j)\pr[Z_{j1} \ge 1]\pr[T_{bn}\uj =l-j]}\non\\
&\le& \sjlt r_j^{-1} j^{-1-2q}\l^2(j) p_0^{-1}KL\l(l)(2/l)^{1+q}\non\\
&&\qquad\qquad\mbox{}
    + \{r^*(l/2)\}^{-1}p_0^{-1}\{L\l(l)\}^2(2/l)^{1+2q} \non\\
&:=& \l(l)l^{-q}\h_3(l), \label{2nd-line-term2.1}
\ena
where clearly $\h_3(l)  = o(1)$ as $l\to\infty$. Then we also have
\eqa
\lefteqn{\sjl  j^{-q}\l(j)\pr[Z_{j1}=1]\pr[T_{bn}\uj =l-2j]}\non\\
&\le& \sjlf r_j^{-1} j^{-1-2q}\l^2(j) p_0^{-1}KL\l(l)(2/l)^{1+q}\non\\
&&\qquad\qquad\mbox{}
    + \{r^*(l/4)\}^{-1} p_0^{-1}\{L^2\l(l)\}^2(4/l)^{1+2q} \non\\
&:=& \l(l)l^{-q}\h_4(l), \label{2nd-line-term2.2}
\ena
again by Lemma~\ref{upper-bnd},
where also $\h_4(l)  = o(1)$ as $l\to\infty$.  The remaining piece
of the last term in~\Ref{2nd-line-term2} is
split into two, as in the proof of the previous lemma, though the
argument is a little simpler.  The bound
\eqa
\lefteqn{\sjl  j^{-q}\l(j) \sum_{s\ge2} \bone_{\{j(s+1) \le \lot\}}
    \pr[Z_{j1}=s]\pr[T_{bn}\uj =l-j(s+1)]}\non\\
&\le& p_0^{-1}KL\l(l)(2/l)^{1+q} \sjl r_j^{-1} j^{-1-2q}\l^2(j)G \e(j)
       \non\\
&:=&  \l(l)l^{-q}\h_5(l), \label{2nd-line-term2.3}
\ena
with $\h_5(l)  = o(1)$ as $l\to\infty$, follows immediately.  For
the second part, we have
\eqa
\lefteqn{\sjl  j^{-q}\l(j) \sum_{s\ge2} \bone_{\{\lot < j(s+1) \le l\}}
    \pr[Z_{j1}=s]\pr[T_{bn}\uj =l-j(s+1)]}\non\\
&\le& p_0^{-1}\sum_{s=2}^{l-1}
   \sum_{j=\lfloor l/2(s+1) \rfloor + 1}^{\lfloor l/(s+1) \rfloor}
   r_j^{-1}\l^2(j) j^{-1-2q} \e(j)\g_s  \pr[T_{bn} =l-j(s+1)]\non\\
&\le& p_0^{-1}\sum_{s=2}^{l-1} \{r^*(l/2(s+1))\}^{-1}
   \e^*(\lfloor l/2(s+1) \rfloor) L^2L_s\l(l) \L_{q/2}
          \{2(s+1)/l\}^{1+3q/2}\g_s\non\\
&\le& \{r^*(0)\}^{-1}\e^*(0)p_0^{-1}3^{1+3q/2}\sum_{s=2}^{l-1}
              L^2L_s\l(l) \L_{q/2}(s/l)^{1+q}\g_s\non\\
&\le& \l(l)l^{-q}\h_6(l),  \label{2nd-line-term2.4}
\ena
with
$$
\h_6(l) :=  \{r^*(0)\}^{-1}\e^*(0)p_0^{-1}3^{1+3q/2}L^2\L_{q/2}G_q l^{-1}
    = o(1)\quad {\rm as}\ l\to \infty.
$$

Combining the results from \Ref{last-line} --~\Ref{2nd-line-term2.4},
it follows from~\Ref{general-probs} that, for $l\le n$,
$$
l\pr[T_{bn}=l] = \sjbl j^{-q}\l(j)\pr[T_{bn} = l-j] + \l(l)l^{-q}\h_7(l),
$$
where $\h_7(l) = o(1)$ as $l\to\infty$.  Hence we deduce that
\eqa
\lefteqn{\l^{-1}(l)l^{1+q} \pr[T_{bn}=l]}\label{asymp.1}\\
&=& \pr[T_{bn} \le l-b-1] +
   \sum_{s=0}^{l-b-1} \Blb {l^q\l(l-s) \over (l-s)^q\l(l)} - 1 \Brb
        \pr[T_{bn} = s] + \h_7(l).\non
\ena
In view of \Ref{lambda-A3}, we can find a sequence~$s_l \to \infty$
such that $s_l = o(l)$ and
$$
\max_{1\le s\le s_l}\left|\frac{\l(l-s)}{\l(l)} - 1\right|
        = o(1) \quad {\rm as}\ l\to \infty:
$$
hence also
$$
\sum_{s=0}^{s_l} \Blm {l^q\l(l-s) \over (l-s)^q\l(l)} - 1 \Brm
        \pr[T_{bn} = s] = \h_8(l) = o(1)  \quad {\rm as}\ l\to \infty.
$$
It then follows from~\Ref{lambda-A2} and~\Ref{T0inf-finite} that
\eqa
&& \sum_{s=s_l+1}^{\lot} \Blm {l^q\l(l-s) \over (l-s)^q\l(l)} - 1 \Brm
        \pr[T_{bn} = s]
               \le (2^qL+1)\pr[T_{bn} > s_l] \non\\
&&\qquad\qquad\le (2^qL+1)\pr[T_{0\infty} > s_l]
= \h_{9}(l) = o(1)
     \quad {\rm as}\ l\to \infty.\non
\ena
For the remaining sum, we use Lemma~\ref{upper-bnd} to give
\eqa
&&\sum_{s=\lot+1}^{l-b-1}
  \Blm {l^q\l(l-s) \over (l-s)^q\l(l)} - 1\Brm \pr[T_{bn} = s]
    \non\\
&&\qquad\le KL\l(l) (2/l)^{1+q}\Blb {l\over2}
            + \sum_{s=1}^{\lot} {l^q\l(s) \over \l(l) s^q}\Brb \non\\
&&\qquad\le KL2^q\L_{q/2}\Blb l^{-q/2} + (2/l)\sum_{s=1}^{\lot}s^{-q/2}\Brb\non\\
&&\qquad = \h_{10}(l) = o(1) \quad{\rm as}\ l \to \infty.
\ena

Putting these estimates into~\Ref{asymp.1}, it follows
that, for $1\le l\le n$,
\eq\label{eps-11}
\l^{-1}(l)l^{1+q} \pr[T_{bn}=l] = 1 - \pr[T_{bn} > l-b-1] + \h_{11}(l),
\en
where $\h_{11}(l) = o(1)$
as $l\to\infty$.  Finally, since also, for $b \le \lot$,
$$
\pr[T_{bn} > l-b-1] \le \pr[T_{0\infty} > l/2] \to 0
    \quad {\rm as}\ l\to\infty,
$$
whereas, for $\lot < b < l$,
\eqa
\pr[T_{bn} > l-b-1] &\le& \pr[T_{b\infty} > 0]
                      \le \pr[ T_{\lot,\infty} > 0]\non\\
&\le& \sum_{j=\lot}^\infty \l(j)j^{-1-q} \to 0
    \quad {\rm as}\ l\to\infty,\non\\
\ena
it follows from~\Ref{eps-11} that, for all $n\ge l$ and $0 \le b \le l-1$,
we have
$$
|\l^{-1}(l)l^{1+q}\pr[T_{bn}=l] - 1| \le H(l),
$$
where $\lim_{l\to\infty}H(l) = 0$, as required.
\epr

\begin{remark}\label{p0-rk}
\textnormal{
The assumption~\Ref{P0-bnd}, that $p_0 > 0$,
can be dispensed with, whatever the distributions of the~$Z_j$,
provided that~\Ref{epsilon-A2} holds.  Clearly, for some $m\ge1$
and $t_1,\ldots,t_m$,
we have
$$
p_0' := \min\Blb \min_{j\ge m+1}\pr[Z_{j1}=0],
       \min_{1\le j\le m}\pr[Z_{j1}=t_j] \Brb > 0,
$$
since $\lim_{j\to\infty} \ex Z_j = 0$.
Then, for $j\le m$ and $s > t_j$, we have the simple bound
$$
\pr[T_{bn}\uj = l-js] \le \pr[T_{bn} = l-j(s-t_j)]/\pr[Z_{j1}=t_j],
$$
which can be used as before, together with the induction
hypothesis, to bound the right hand side of~\Ref{general-probs}
in the proof of Lemma~\ref{upper-bnd}, provided that $s > t_j$.
So, recalling~\Ref{star} with $g = \bone_{\{l\}}$, we write
\eqa
\lefteqn{\sum_{j=1}^m \ex\{jZ_j\bone_{\{l\}}(T_{bn})\}} \non\\
&=& \sum_{j=1}^m jr_j
       \ex\{Z_{j1}I[Z_j \le t_j]\bone_{\{l\}}(T_{bn})\}\non\\
&&\qquad \mbox{} + \sum_{j=1}^m jr_j
      \ex\{Z_{j1}I[Z_j > t_j]\bone_{\{l\}}(T_{bn})\}.\non
\ena
The second term is estimated exactly as before.  The first is
no larger than $\kappa\pr[T_{bn}=l]$, where
$$
\kappa := \sum_{j=1}^m jr_jt_j,
$$
and hence can be taken onto the left hand side of~\Ref{first-bnd}
whenever $l \ge 2\kappa$; with these modifications, the proof
of Lemma~\ref{upper-bnd} can be carried through as before.
The proof of Theorem~\ref{LLT} requires almost no modification,
if~$p_0$ is replaced by~$p_0'$.}
\end{remark}

\begin{remark}\label{id-rk}
\textnormal{
If the $(Z_j,\,j\ge1)$ are infinitely divisible,
then we can choose the~$r_j$ to be arbitrarily large for each
fixed~$j$, in the limit
making $\h_k(l) = 0$, $3 \le k \le 6$, and $p_0 = 1$.  The limiting
values as $r_j \to \infty$ of $\e_{js}$, for fixed~$j$ and~$s\ge1$,
are {\it not\/} however in general zero.}
\end{remark}

\begin{remark}\label{Gq-rk}
\textnormal{
  The assumption~\Ref{epsilon-A3} that~$G_q$
be finite is not just an artefact of the proofs.
It appears in particular when bounding
the quantity~$S_2$ in~\Ref{S2-bnd} in the proof of
Lemma~\ref{upper-bnd}, and is an element in the
quantity~$\h_2(l)$, which contributes to the bound on~$H(l)$ in
Theorem~\ref{LLT}.  However, $l^{-1}S_2$ is of the same order as the
probability that~$T_{0n}$ is composed of~$s$ components of
equal sizes~$\lfloor l/s \rfloor$, plus a small remainder,
for some~$s\ge2$, and $G_q < \infty$ is the condition which
ensures that this probability is of smaller order than
$\l(l)l^{-1-q}$.}
\end{remark}

\subsection{Proof of Theorem~\ref{joint-dns}}
As in [ABT, Lemma~3.1], it follows from the Conditioning Relation that,
for any $b \le n$,
\eqa
&&d_{TV}(\law(C_1\un,\ldots,C_{b}\un),\law(Z_1,\ldots,Z_b))\non\\
&&\qquad = \sum_{j\ge0} \pr[T_{0b}=j]
    \Blb 1 - {\pr[T_{bn}=n-j] \over \pr[T_{0n}=n]} \Brb_+. \label{TV}
\ena
Pick $b=b(n)$ with $n-b(n) \to \infty$, and observe that the right hand side
of~\Ref{TV} is at most
$$
\pr[T_{0b} > j_n] + \ex g_n(T_{0b}),
$$
where $g_n(j) = 0$ for $j > j_n$ and where, for all~$n$ such that
$H(n) < 1/2$,
$$
0 \le g_n(j) \le \Blm {n^{1+q}\l(n-j) \over (n-j)^{1+q}\l(n)} - 1 \Brm
    + 2^{1+q}L\,2(H(n) + H(n-j)), \quad 0 \le j \le j_n,
$$
from Theorem~\ref{LLT},
provided that $0 \le j_n \le \lfloor n/2 \rfloor$ and that $j_n \le
n-b(n) -1$. This implies in particular that~$g_n(j)$ is
uniformly bounded for sequences~$j_n$ satisfying these conditions.
Now, from \Ref{lambda-A3} and Theorem~\ref{LLT}, it follows that
$\lim_{n\to\infty} g_n(j) = 0$ for each fixed~$j$.
Since also $T_{0b} \le T_{0\infty}$ a.s.\ and $T_{0\infty}$ is
a.s.\ finite, it follows by dominated convergence that
$\lim_{n\to\infty}\ex g_n(T_{0b(n)}) = 0$, provided that,
in the definition of~$g_n$,
$j_n \le \min\{n-b(n) - 1, \lfloor n/2 \rfloor\}$.  On the
other hand,
$$
\pr[T_{0b(n)} > j_n] \le \pr[T_{0\infty} > j_n] \to 0,
$$
so long as $j_n \to \infty$.  Thus, taking for example
$b(n) = \lfloor 3n/4 \rfloor$ and $j_n = \lfloor n/4 \rfloor - 1$,
it follows that
$$
d_{TV}(\law(C_1\un,\ldots,C_{b(n)}\un),\law(Z_1,\ldots,Z_{b(n)}))
  \to 0
$$
as $n\to\infty$.  On the other hand,
we have $\sum_{j= \lfloor 3n/4 \rfloor +1}^n C_j\un \le 1$ a.s.,
because $T_{0n}(C\un) = n$ a.s., by the definition of $C\un$.
Hence, with~$b(n)$ as above, we have $C_j\un = 0$ a.s.\ for all
$j> b(n)$ if $T_{0b(n)}(C\un) = n$, while if
$T_{0b(n)}(C\un) = t$ for some $t
< n-b(n)$, then $C_{n-t}\un = 1$
and $C_j\un = 0$ for all other $j>b(n)$.  This proves the theorem.
\epr

\bsk\bsk
\nin{\bf Acknowledgements}

\msk
ADB was supported in this research by
Schweizer Nationalfondsprojekt Nr.~20-67909.02.
BG's research was supported by the Fund for the
Promotion of Research at Technion.  We wish to thank the referees
for a number of helpful suggestions.

\bsk\bsk




\end{document}